\documentclass{snapshotmfo}

\junioreditor[gs]{will be filled out by the editors}{junior-editors@mfo.de}
\senioreditor[f]{Anja Randecker}{senior-editor@mfo.de}
\director[m]{Gerhard Huisken}

\usepackage[utf8]{inputenc}
\usepackage{amsmath,amssymb}
\usepackage{tikz}
\usetikzlibrary{arrows.meta}
\usetikzlibrary{calc}
\usepackage[all]{xy}
\usepackage{wasysym}
\usepackage{xcolor}

\usepackage[USenglish]{babel}

\usepackage[noabbrev,capitalise,nameinlink]{cleveref} 

\usepackage{ellipsis}

\newenvironment{statement}{\begin{quote}}{\end{quote}}

\author{Darij Grinberg \and Tom Roby}

\title{Detropicalization as a proof technique}

\authorinfo{\authorname{Darij Grinberg} is a professor of mathematics at Drexel University.}
\authorinfo{\authorname{Tom Roby} is a professor of mathematics at University of Connecticut.}

\begin{document}
\pdfbookmark[-1]{Detropicalization as a proof technique}{snapshottitle}

\begin{abstract}
Rational functions make sense over any commutative ring, as long as the denominators are invertible. When there is no subtraction involved, they even apply over semirings (rings without subtraction). It is particularly worthwhile to evaluate them over the \emph{tropical semiring}, in which the roles of addition and multiplication are played by maxima and addition. Over this semiring, algebraic results often acquire combinatorial meaning. We give a few examples.
\end{abstract}


\section{Introduction}

Tropical algebra begins with a simple observation about the laws
of arithmetic:
The distributive law $\left(  a+b\right)  c=ac+bc$ for
multiplication over addition resembles the distributive law $\max\left\{
a,b\right\}  +c=\max\left\{  a+c,b+c\right\}  $ for addition over
the maximum function.
This is just one of several instances of a parallelism between the
pair of operations $\left(  \cdot,+\right)  $
and the pair of operations $\left(  +,\max\right)  $.
More such similarities are easily found, although there are limits
to this parallelism (e.g., addition can be undone by subtraction,
but taking a maximum cannot be undone).

As an algebraist, one is thus tempted to axiomatize the common
parts of the theories of $\left(  \cdot,+\right)  $ and
$\left(  +,\max\right)  $, and use this new abstraction as a
bridge to transfer results and intuitions from one theory to the other.

This is the chief idea behind the notions of
\emph{tropicalization} and \emph{detropicalization} and the modifier
\emph{tropical} attached to several mathematical disciplines (most commonly
\emph{tropical geometry} \cite{Brugalle-ss}). The abstraction that unites the
$\left(  \cdot,+\right)  $ and $\left(  +,\max\right)  $ settings is the
concept of a \emph{semiring}. A \emph{semiring} $R$ is defined just like a
ring, but without requiring a subtraction (so its addition merely makes $R$ a
monoid, not a group); thus, examples of semirings include $\mathbb{N}=\left\{
\text{nonnegative integers}\right\}  $, $\mathbb{R}_{\geq0}=\left\{
\text{nonnegative reals}\right\}$,  as well as many versions of nonnegative
polynomials (e.g., polynomials with nonnegative coefficients; polynomials
taking only nonnegative values on $\mathbb{R}_{\geq0}$ or on $\mathbb{R}$;
polynomials writable as sums of squares, etc.). Of course, any ring is also
a semiring. The simplest example of a (non-ring) semiring is the
\emph{Boolean semiring} $\mathbb{B}=\left\{  0,1\right\}  $ with $1+1=1$
(!) but otherwise the usual rules of arithmetic. As its name suggests, its
elements $0$ and $1$ model the truth values \textquotedblleft
False\textquotedblright\ and \textquotedblleft True\textquotedblright, with
$+$ modelling \textquotedblleft or\textquotedblright\ and $\cdot$ modelling
\textquotedblleft and\textquotedblright.

Just a bit deeper lies the \emph{(max-plus) tropical semiring}
$\mathbb{T}_{\mathbb{R}}$ (of real numbers), which as a set is $\mathbb{R}%
\cup\left\{  -\infty\right\}  $ (that is, the real numbers with a symbol
$-\infty$ thrown in) but has the operations $\left(  x,y\right)  \mapsto
\max\left\{  x,y\right\}  $ and $\left(  x,y\right)  \mapsto x+y$ acting as
``addition''\ and
``multiplication''\ (where $\max\left\{-\infty, x\right\} := x$ and
$\left(-\infty\right) + x := -\infty$ for all $x$).
Thus, the new element $-\infty$ takes the
role of \textquotedblleft zero\textquotedblright, while $0\in\mathbb{R}$
becomes \textquotedblleft one\textquotedblright. This structure is
a commutative semiring (e.g., we saw above that it is distributive).
Even better,
the tropical semiring $\mathbb{T}_{\mathbb{R}}$ also has a \textquotedblleft
division\textquotedblright\ (the operation $\left(  x,y\right)  \mapsto x-y$,
valid whenever $y$ is \textquotedblleft nonzero\textquotedblright, i.e.,
distinct from $-\infty$), and hence becomes a \emph{semifield}. The same
construction can be done with any ordered abelian group instead of
$\mathbb{R}$, such as $\mathbb{R}^{k}$ with lexicographic order; the operation
$\max$ can also be replaced by $\min$ (with $-\infty$ replaced by $+\infty$).

Computer scientists have used tropical semirings to study formal languages
\cite{Pin-tropsemi} and networks \cite{Maxplus}. Algebraic geometers have
taken an interest in \emph{tropical geometry}, which models some qualitative
aspects of algebraic varieties using their analogues over
$\mathbb{T}_{\mathbb{R}}$; see, e.g., \cite{SpSt-TM}, \cite{Brugalle-sa},
\cite{Brugalle-ss}, \cite{MaclaganSturmfels}, \cite{GKMV}.
The former often proceed from
the tropical setting and generalize to semirings for clarity, whereas the
latter start with objects defined over rings and fields and transfer them to
the tropical world via the general realm of semirings. Some of the resulting
tropical concepts have turned out to be known under different names, such as
phylogenetic trees \cite{SpSt-TM}, \cite[\S 2]{RiTrYu-NAMS},
\cite{PachterSturmfels}.

More recently, the same type of knowledge transfer has caught on in algebraic
combinatorics. While continuous piecewise-linear functions (the tropical
analogue of rational functions) have been known to combinatorialists for a
long time, it was probably Kirillov \cite{Kirillov-ITC} (in joint work with
Arkadiy Berenstein, inspired by the work of Lusztig and Kashiwara on
combinatorial crystals) who started viewing them as tropical rational
functions and generalizing them to arbitrary semirings -- a process
he dubbed \emph{tropification} but is now better known as
\emph{detropicalization}. It is this process that we shall illustrate and
give examples for in this snapshot.

\section{The racetrack puzzle}

We begin with a simple yet surprising example. We recall a well-known puzzle
(e.g., \cite[Honsberger, \textquotedblleft The track\textquotedblright, p.
283]{Klarner-Gardener}):

\begin{statement}
\textbf{The racetrack puzzle:} A circular track has $n$ gas cans on it. A car
with an empty tank needs to make a single lap (clockwise) around the track. Taken together,
the $n$ cans contain just enough gas for the car 
to make such a lap. Prove that there is at least one starting point on the
track such that, if the car is airlifted to this point, then (by refueling
along the way) it can make such a lap.\footnote{The car's tank is large enough
to fit all the gas.}
\end{statement}

\begin{figure}[ht]
\centering
\begin{tikzpicture}[scale=2.2]
\node(a) at ({cos(0)}, {sin(0)}) {};
\node(b) at ({cos(-0.1*360)}, {sin(-0.1*360)}) {};
\node(c) at ({cos(-0.4*360)}, {sin(-0.4*360)}) {};
\node(d) at ({cos(-0.5*360)}, {sin(-0.5*360)}) {};
\node(e) at ({cos(-0.8*360)}, {sin(-0.8*360)}) {};
\fill(a) circle (1pt) node[right] {$A$} node[left]{$(0.2)$};
\fill(b) circle (1pt) node[below right] {$B$} node[above left]{$(0.1)$};
\fill(c) circle (1pt) node[below left] {$C$} node[above right] {$(0.4)$};
\fill(d) circle (1pt) node[left] {$D$} node[right] {$(0.2)$};
\fill(e) circle (1pt) node[above] {$E$} node[below] {$(0.1)$};
\draw[black] (1,0)
arc[radius = 1, start angle=0, delta angle=0.1*360]
node[above right] {$0.2$}
arc[radius = 1, start angle=0.1*360, end angle=0.35*360]
node[above left] {$0.3$}
arc[radius = 1, start angle=0.35*360, end angle=0.55*360]
node[left] {$0.1$}
arc[radius = 1, start angle=0.55*360, end angle=0.75*360]
node[below] {$0.3$}
arc[radius = 1, start angle=0.75*360, end angle=0.95*360]
node[right] {$0.1$}
arc[radius = 1, start angle=0.95*360, end angle=360];
\draw[red, <<-, very thick] (2,0) arc[radius=2, start angle=0, end angle=20];
\end{tikzpicture}
\caption{Example for the racetrack puzzle ($n = 5$). The five gas cans
$A,B,C,D,E$ have respectively $0.2,0.1,0.4,0.2,0.1$ gallons of gas (numbers in
parentheses). The numbers on the arcs show how much gas is needed to go from
one can to the next.}
\label{fig.track1}
\end{figure}
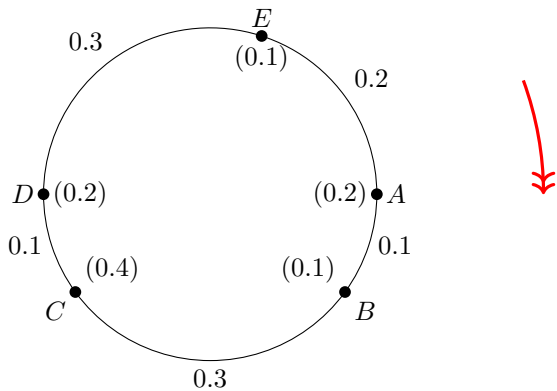

An example is shown on Fig. \ref{fig.track1}. Here, a car starting at gas can
$A$ would reach $B$ with an empty tank, refuel to $0.1$ gas,
but then would get stuck on the way
to $C$, since the total amount of gas in $A$ and $B$ does not suffice for the
trek from $A$ to $C$. On the other hand, a car starting at $C$ would make it
all the way around the track.

Now to something (seemingly) completely different: an algebra problem, which
also has appeared in various contests\footnote{One of us fondly remembers it
from his mathematical olympiad days, gratuitously obfuscated by being stated
as an inequality. We have not managed to find the original source, but a
variant has appeared in the Final Round of the Swiss Mathematical Olympiad
2005 \cite[Finalrunde 2005, Aufgabe 6]{SMO-AlPr}.}:

\begin{statement}
\textbf{The three fractions puzzle:} For any three real numbers $a,b,c$ satisfying
$abc=1$, prove that%
\[
\dfrac{1}{1+a+ab}+\dfrac{1}{1+b+bc}+\dfrac{1}{1+c+ca}=1
\]
(provided that all denominators are nonzero).
\end{statement}

Of course, in this form, this is not exactly challenging, as one can simply
cross-multiply and expand. However, the more elegant solutions have the
advantage of generalizing from three to $n$ numbers, thus solving the
following generalization of the puzzle:

\begin{statement}
\textbf{The $n$ fractions puzzle:} Given $n$ real numbers $a_{1},a_{2}%
,\ldots,a_{n}$ satisfying $a_{1}a_{2}\cdots a_{n}=1$, prove that%
\begin{equation}
\sum_{i=1}^{n}\dfrac{1}{1+a_{i}+a_{i}a_{i+1}+\cdots+a_{i}a_{i+1}\cdots
a_{i+n-2}}=1, \label{eq.nfr}%
\end{equation}
where the indices are cyclic modulo $n$ (so $a_{j+n}=a_{j}$ for each $j$),
and where we assume that all denominators are nonzero.
\end{statement}

Even in this form, the puzzle is not too difficult (see the Appendix for a
solution). But it is far less well-known that it is closely related to
the racetrack puzzle. How?

It helps to restate the racetrack puzzle. We denote the $n$ gas cans by
$C_{1},C_{2},\ldots,C_{n}$, the amount of gas available at can $C_{i}$ by
$c_{i}$, and the amount of gas necessary to get from $C_{i}$ to $C_{i+1}$ by
$d_{i}$ (where the indices are cyclic modulo $n$). Then, a car that starts at
the can $C_{i}$ will

\begin{itemize}
\item reach the can $C_{i+1}$ if $c_{i}\geq d_{i}$;

\item then reach the can $C_{i+2}$ if $c_{i}+c_{i+1}\geq d_{i}+d_{i+1}$;

\item then reach the can $C_{i+3}$ if $c_{i}+c_{i+1}+c_{i+2}\geq
d_{i}+d_{i+1}+d_{i+2}$;

\item and so on.
\end{itemize}

Hence, this car will finish its lap if and only if%
\[
c_{i}+c_{i+1}+\cdots+c_{i+k-1}\geq d_{i}+d_{i+1}+\cdots+d_{i+k-1}%
\qquad\text{for each }1\leq k\leq n.
\]
Setting $a_{j}:=d_{j}-c_{j}$, we can rewrite this condition as%
\begin{equation}
a_{i}+a_{i+1}+\cdots+a_{i+k-1}\leq0\qquad\text{for each }1\leq k\leq n.
\label{eq.cond.3}%
\end{equation}

On the other hand, we are given that the $n$ cans have just enough fuel for
the entire track; in our notations, this means that%
\[
c_{1}+c_{2}+\cdots+c_{n}=d_{1}+d_{2}+\cdots+d_{n},
\]
or, equivalently,%
\begin{equation}
a_{1}+a_{2}+\cdots+a_{n}=0. \label{eq.cond.0}%
\end{equation}
Thus, by cyclic periodicity, $a_{i}+a_{i+1}+\cdots+a_{i+n-1}=0$ for all $i$.
Hence, in the condition (\ref{eq.cond.3}), we can replace the $k=n$ case by
the $k=0$ case (in which case $a_{i}+a_{i+1}+\cdots+a_{i+k-1}$ is an empty sum
and equals $0$ by definition). Hence, we can rewrite the condition
(\ref{eq.cond.3}) as%
\begin{equation}
a_{i}+a_{i+1}+\cdots+a_{i+k-1}\leq0\qquad\text{for each }0\leq k\leq n-1,
\label{eq.cond.4}%
\end{equation}
i.e., as%
\begin{equation}
\max\nolimits_{0\leq k\leq n-1}\left\{  a_{i}+a_{i+1}+\cdots+a_{i+k-1}%
\right\}  =0 \label{eq.cond.5}
\end{equation}
(since this maximum cannot be smaller than $0$, seeing that $0$
belongs to the set).

Thus we aim to show that there exists some gas can
$C_i$ that satisfies
(\ref{eq.cond.5}). Equivalently, we aim to show that%
\footnote{Why ``equivalently''?
Well, proving (\ref{eq.cond.6}) will clearly yield (\ref{eq.cond.5}).
But the converse implication is also easy, since the minimum
in \eqref{eq.cond.6} cannot be smaller than $0$.}
\begin{equation}
\min\nolimits_{1\leq i\leq n}\left\{  \max\nolimits_{0\leq k\leq n-1}\left\{
a_{i}+a_{i+1}+\cdots+a_{i+k-1}\right\}  \right\}  =0.
\label{eq.cond.6}
\end{equation}
We can rewrite the minimum here in terms of maxima by recalling that
we have
$\min\left\{  u_{1},u_{2},\ldots,u_{k}\right\}  =-\max\left\{  -u_{1}%
,-u_{2},\ldots,-u_{k}\right\}  $ for any real
numbers $u_{1},u_{2},\ldots,u_{k}$;
thus, (\ref{eq.cond.6}) takes the equivalent form%
\[
-\max\nolimits_{1\leq i\leq n}\left\{  -\max\nolimits_{0\leq k\leq
n-1}\left\{  a_{i}+a_{i+1}+\cdots+a_{i+k-1}\right\}  \right\}  =0,
\]
i.e., the equivalent form%
\begin{equation}
\max\nolimits_{1\leq i\leq n}\left\{  -\max\nolimits_{0\leq k\leq n-1}\left\{
a_{i}+a_{i+1}+\cdots+a_{i+k-1}\right\}  \right\}  =0. \label{eq.cond.7}%
\end{equation}
So we aim to prove the equality (\ref{eq.cond.7}) under the assumption of
(\ref{eq.cond.0}).

Now, recall that the operations $+$ and $\max$ of the ordered abelian group
$\mathbb{R}$ are the multiplication and the addition of the tropical semiring
$\mathbb{T}_{\mathbb{R}}$. Moreover, the operation $x\mapsto-x$ of
$\mathbb{R}$ corresponds to taking reciprocals in $\mathbb{T}_{\mathbb{R}}$.
Thus, in terms of the semiring $\mathbb{T}_{\mathbb{R}}$, our claim
(\ref{eq.cond.7}) is an instance of the identity%
\begin{equation}
\sum_{i=1}^{n}\dfrac{1}{\sum_{k=0}^{n-1}a_{i}a_{i+1}\cdots a_{i+k-1}}=1,
\label{eq.cond.7t}%
\end{equation}
whereas our assumption (\ref{eq.cond.0}) is saying that $a_{1}a_{2}\cdots
a_{n}=1$. But (\ref{eq.cond.7t}) is precisely (\ref{eq.nfr}). Hence, the
racetrack puzzle is a particular case of the $n$ fractions puzzle, just
with real numbers replaced by the elements of an arbitrary commutative
semiring
(assuming that the denominators are invertible)! (The proof of the $n$
fractions puzzle we give in the Appendix works in this generality.)

\section{Detropicalization}

The above example illustrates the technique of \emph{detropicalization}:
proving a combinatorial result (in our case, the racetrack puzzle) by
rewriting it as an identity in the tropical semiring (a \textquotedblleft
max-plus identity\textquotedblright) and then algebraically proving this
identity for \textbf{all} semirings, or, more commonly, all commutative
semirings.\footnote{Actually it suffices to prove it for the ring
$\mathbb{R}$; its validity for the tropical semiring
$\mathbb{T}_\mathbb{R}$ then follows by a limiting process.
See \cite[\S 7]{EinsteinPropp-paper}
for details (note that $\min$ is used instead of $\max$ in
\cite{EinsteinPropp-paper}, but this is clearly isomorphic via
$x \mapsto -x$).}

This technique is far from a panacea. For one, not every combinatorial claim
can be restated as a max-plus identity. Next, such an identity might not
remain valid when generalized from the tropical semiring to arbitrary
semirings\footnote{In fact, \textbf{most} identities don't! The identity
$\left(  x+y\right)  ^{2}=x^{2}+y^{2}$ holds for every tropical semiring, but
fails for most other semirings. For a more combinatorial
example, consider the map $\operatorname*{sort}:\mathbb{R}^{2}\rightarrow
\mathbb{R}^{2}$ that sends each pair $\left(  x,y\right)  $ to $\left(
\min\left\{  x,y\right\}  ,\max\left\{  x,y\right\}  \right)  $. In the
combinatorial world, this map is idempotent (i.e., $\operatorname*{sort}%
\circ\operatorname*{sort}=\operatorname*{sort}$), but the obvious
generalization to semirings (replacing $\max\left\{  x,y\right\}  $ by $x+y$
and $\min\left\{  x,y\right\}  $ by $\dfrac{1}{\dfrac{1}{x}+\dfrac{1}{y}}$)
has no such property.}. Finally, the generalization won't always be easier to
prove than the original claim. However, more often than not, if such a
generalization is true, then its proof won't be difficult, since algebra is
almost always easier than combinatorics. Thus, in our experience, an attempt
at detropicalization that makes it through the first two gates (restatement
and correctness) will make through the third as well.

At the origins of detropicalization stands the striking success of Kirillov
\cite{Kirillov-ITC} (see also \cite{NoumiYamada})
detropicalizing several combinatorial operations on
semistandard tableaux: the RSK algorithm, Sch\"{u}tzenberger promotion and
evacuation and the Lascoux--Sch\"{u}tzenberger action. These are nowadays
regarded as parts of the theory of crystal bases \cite{BumpSchilling},
which owes much of its inspiration to the original work by
Lusztig and Kashiwara in the early 1990s that anticipates the idea of
detropicalization to an extent.

We have had two successes applying the detropicalization technique in our
research. One was detropicalizing the properties of \emph{piecewise-linear
rowmotion} on finite posets \cite{BR1}, \cite{BR2} (which we generalized even
further to noncommutative semirings in \cite{NCBR}); another was the proof of
the Pelletier--Ressayre identities for Littlewood--Richardson coefficients
\cite{Grinberg-PR}. In both cases, we had a map defined on combinatorial
objects, and were trying to show certain identities for this map (e.g., that a
certain power of this map equals the identity). The map could be defined in
terms of maxima and addition of integers, i.e., it was the tropical case of a
universal map (natural transformation) defined for any commutative semiring;
and the identities in question turned out to hold in that generality.
Moreover, proving them in this generality turned out to be easier than doing
so in the tropical semiring alone.

This last fact might be surprising, so here is some elaboration. First, as
algebraists, we were simply more used to working with sums and products than
with maxima and sums, and the very fact of moving from combinatorics to
algebra gave us the feeling of treading on stable ground. But beside this, the
category of commutative semirings has several more substantial advantages,
including the existence of free objects; in down-to-earth terms, this entails
that an identity holds for all commutative semirings if it holds for the
polynomial rings $\mathbb{N}\left[  x_{1},x_{2},\ldots,x_{k}\right]  $ of
multivariate polynomials with coefficients in $\mathbb{N}$. In particular, in
proving an identity for all commutative semirings, we could use subtraction,
since $\mathbb{N}\left[  x_{1},x_{2},\ldots,x_{k}\right]  $ embeds into
the ring $\mathbb{Z}\left[  x_{1},x_{2},\ldots,x_{k}\right]  $. This allowed
us to use determinants and Zariski density arguments in \cite{BR1}, \cite{BR2},
and to use an inequality trick in \cite[Remark 3.16]{Grinberg-PR}. None of
these techniques would have been available in the original tropical setting.

We had less luck detropicalizing the Ford--Fulkerson
algorithm, despite it being describable as a tropical rational function
\cite{MO-FF}; perhaps semirings are not its natural habitat.
However, detropicalization is known to be useful in graph theory; e.g.,
taking powers of the adjacency matrix over the tropical semiring leads to a
shortest-path algorithm \cite[\S 3.3]{KepnerGilbert}.

\section{Somos-4 and the Laurent phenomenon}

We can also work from the opposite direction: Find a property of commutative
semirings and apply it to the tropical semiring; reinterpret the result
combinatorially. Mathematically, nothing can go wrong here; of course,
the resulting interpretation may or may not be interesting.
We give a brief example.%
\footnote{We were informed that this example had also been found by Sergey
Fomin and Pavlo Pylyavskyy.}

Let $F$ be the field of rational functions in six variables $p,q,r,s,w,z$ over
$\mathbb{Q}$. Define a sequence $\left(  a_{0},a_{1},a_{2},\ldots\right)  $ of
elements of $F$ (known as the \emph{Somos-4 sequence with parameters})
recursively by%
\begin{align*}
a_{0} &  =p,\ \ \ \ \ \ \ \ \ \ a_{1}=q,\ \ \ \ \ \ \ \ \ \ a_{2}%
=r,\ \ \ \ \ \ \ \ \ \ a_{3}=s,\\
a_{n} &  =\dfrac{wa_{n-1}a_{n-3}+za_{n-2}^{2}}{a_{n-4}}%
\ \ \ \ \ \ \ \ \ \ \text{for }n\geq4.
\end{align*}
One instance of the famous \emph{Laurent phenomenon} (see \cite[Example
10.5]{FordyMarsh}) says that each $a_{n}$ is a Laurent polynomial in
$p,q,r,s,w,z$ (that is, a polynomial in $p,q,r,s,w,z$ divided by a monomial).
From the positivity results of Lee--Schiffler \cite[Theorem 1.1]{LeeSchiffler}
and Gross--Hacking--Keel--Kontsevich \cite[Corollary 0.4]{GHKK}, we further
know that the coefficients of this Laurent polynomial are nonnegative integers.

Since this holds in $F$, this will also hold in any commutative semiring,
thus in particular in the tropical semiring $\mathbb{T}_{\mathbb{R}}$.
In other words, if $p,q,r,s,w,z$ are real numbers, and if we define a
sequence $\left(  a_{0},a_{1},a_{2},\ldots\right)  $ of real numbers
recursively by%
\begin{align*}
a_{0} &  =p,\ \ \ \ \ \ \ \ \ \ a_{1}=q,\ \ \ \ \ \ \ \ \ \ a_{2}%
=r,\ \ \ \ \ \ \ \ \ \ a_{3}=s,\\
a_{n} &= \max\left\{  w+a_{n-1}+a_{n-3},\ z+2a_{n-2}\right\}  -a_{n-4}
\ \ \ \ \ \ \ \ \ \ \text{for }n\geq4,
\end{align*}
then each $a_{n}$ is a \textquotedblleft tropical Laurent
polynomial\textquotedblright\ in $p,q,r,s,w,z$. Of course, a tropical Laurent
polynomial is what comes out when a Laurent polynomial is interpreted in the
tropical semiring -- i.e., a pointwise maximum of a finite set of linear
forms. Note, in particular, that such tropical Laurent polynomials are always
convex in their inputs, so we conclude that $a_{n}$ (in the tropical setting)
is convex as a function of the vector $\left(  p,q,r,s,w,z\right)  $.

\section{Birational rowmotion on rectangles}

We now present a more recent application of detropicalization techniques,
which has become known as \emph{birational rowmotion} although closely related
systems appear in various other guises (birational promotion, some Y- and
T-systems, octahedron recurrences \cite{Roby16}, \cite{JosLiu-oct}).

The simplest way to define it is as a \textquotedblleft discrete dynamical
system\textquotedblright. Consider two reals $a$ and $b$, two positive
integers $p$ and $q$, and an infinite family $\left(  x_{i,j,k}\right)
_{i\in\left[  p\right]  ,\ j\in\left[  q\right]  ,\ k\in\mathbb{Z}}$ of reals,
where $\left[  n\right]  :=\left\{  1,2,\ldots,n\right\}  $ for each
$n\in\mathbb{N}$. We can think of this family as a (doubly) infinite tower of
$p\times q$-matrices $\left(  x_{i,j,k}\right)  _{i\in\left[  p\right]
,\ j\in\left[  q\right]  }$, one for each integer $k$. We say that this family
is a \emph{rowmotion family} if it satisfies%
\begin{align}
& x_{i,j,k+1}+x_{i,j,k}\nonumber\\
& =\max\left\{  x_{i,j-1,k},\ x_{i-1,j,k}\right\}  +\min\left\{
x_{i,j+1,k+1},\ x_{i+1,j,k+1}\right\}  \label{eq.pwr}%
\end{align}
for all $i,j,k$, where we agree that

\begin{itemize}
\item if one of the two elements of the set $\left\{  x_{i,j-1,k}%
,\ x_{i-1,j,k}\right\}  $ is undefined (i.e., one of $i-1$ and $j-1$ is $0$),
then the maximum of this set is the other element;

\item if both are undefined, then the maximum of this set is $a$;

\item if one of the two elements of the set $\left\{  x_{i,j+1,k+1}%
,\ x_{i+1,j,k+1}\right\}  $ is undefined, then the minimum is the other element;

\item if both are undefined, then the minimum is $b$.
\end{itemize}

Equation~(\ref{eq.pwr}) can be solved uniquely by recursion if all the
values $x_{i,j,0}$ (along with $a$ and $b$) are given: Indeed, knowing all
$x_{i,j,k}$ for a given $k$ allows us to compute all $x_{i,j,k+1}$ in the
order of decreasing $i+j$, whereas the $x_{i,j,k-1}$ can be computed in the
order of increasing $i+j$. Thus, (\ref{eq.pwr}) really defines an invertible
piecewise-linear map $\rho_{\operatorname*{PL}}:\mathbb{R}^{p\times
q}\rightarrow\mathbb{R}^{p\times q}$ that sends the matrix $\left(
x_{i,j,k}\right)  _{i\in\left[  p\right]  ,\ j\in\left[  q\right]  }$ to
$\left(  x_{i,j,k+1}\right)  _{i\in\left[  p\right]  ,\ j\in\left[  q\right]
}$ for each $k\in\mathbb{Z}$. This map is known as \emph{piecewise-linear
rowmotion}, and has been introduced by Einstein and Propp
\cite{EinsteinPropp-paper}. They defined it as a sequence of $pq$
\textquotedblleft piecewise-linear toggles\textquotedblright%
\ \cite[Definitions 3.1 and 3.2]{EinsteinPropp-paper}, each of which
corresponds to a single step in solving (\ref{eq.pwr}). They also defined a
similar map called \emph{piecewise-linear promotion} \cite[Definition
3.3]{EinsteinPropp-paper}, which is the map sending $\left(  x_{i,j,k+i}%
\right)  _{i\in\left[  p\right]  ,\ j\in\left[  q\right]  }$ to $\left(
x_{i,j,k+1+i}\right)  _{i\in\left[  p\right]  ,\ j\in\left[  q\right]  }$ for
each $k\in\mathbb{Z}$ (that is, essentially a different \textquotedblleft
slicing\textquotedblright\ of the equation (\ref{eq.pwr})). Both maps
generalize combinatorial constructions: the rowmotion (or
Cameron--Fon-der-Flaass map) on order ideals of a poset, and the promotion of
semistandard Young tableaux. We have here restricted ourselves to families
with rectangle-shaped footprints, but other variants have good behavior as
well (see, e.g., \cite{BR1, BR2}).

Propp and Einstein conjectured\footnote{in the original arXiv version of
\cite{EinsteinPropp-paper} (arXiv:1310.5294). The conjectures had already been proved by 
the time the published version actually appeared some 8 years later.} that both of their maps are periodic
with period $p+q$. In our language, this says that any solution of
(\ref{eq.pwr}) satisfies the \emph{periodicity}%
\[
x_{i,j,k+p+q}=x_{i,j,k}\ \ \ \ \ \ \ \ \ \ \text{for all }i,j,k.
\]
This actually follows from another conjecture: the \emph{reciprocity}
relations%
\[
x_{i,j,k}=a+b-x_{p+1-i,\ q+1-j,\ k-i-j+1}\ \ \ \ \ \ \ \ \ \ \text{for all
}i,j,k.
\]
These conjectures extended a famously nontrivial property of semistandard
tableaux, and appeared even harder.

\begin{figure}[ht]
\centering
\begin{tikzpicture}[
    scale=2,
    vertex/.style={circle, fill=black, inner sep=1.4pt},
    edge/.style={black, line width=0.6pt}
]

\coordinate (u) at (2.2,0.55);
\coordinate (v) at (0.95,-0.45);

\coordinate (h) at (0,1.15);

\def\ext{0.55}

\foreach \k in {0,...,4} {
    \path
        ($\k*(h)$)         node[vertex] (A\k) {}
        ($\k*(h)+(u)$)     node[vertex] (B\k) {}
        ($\k*(h)+(u)+(v)$) node[vertex] (C\k) {}
        ($\k*(h)+(v)$)     node[vertex] (D\k) {};
}

\foreach \k in {0,...,4} {
    \draw[edge] (A\k) -- (B\k) -- (C\k) -- (D\k) -- (A\k);
}

\draw[edge] ($(A0)-\ext*(h)$) -- ($(A4)+\ext*(h)$);
\draw[edge] ($(B0)-\ext*(h)$) -- ($(B4)+\ext*(h)$);
\draw[edge] ($(C0)-\ext*(h)$) -- ($(C4)+\ext*(h)$);
\draw[edge] ($(D0)-\ext*(h)$) -- ($(D4)+\ext*(h)$);

\node[left]        at (A0) {$x$};
\node[above left]        at (B0) {$z$};
\node[right]        at (C0) {$y$};
\node[below left]        at (D0) {$w$};

\node[left]        at (A1) {$\frac{b(x+y)w}{xz}$};
\node[anchor=east, xshift=1pt, yshift=7pt]        at (B1) {$\frac{b(x+y)}{z}$};
\node[right]        at (C1) {$\frac{b(x+y)w}{yz}$};
\node[below left]        at (D1) {$\frac{ab}{z}$};

\node[left]        at (A2) {$\frac{ab}{y}$};
\node[anchor=east, xshift=1pt, yshift=7pt]        at (B2) {$\frac{bw(x+y)}{xy}$};
\node[right]        at (C2) {$\frac{ab}{x}$};
\node[below left]        at (D2) {$\frac{az}{x+y}$};

\node[left]        at (A3) {$\frac{ayz}{w(x+y)}$};
\node[anchor=east, xshift=1pt, yshift=7pt]        at (B3) {$\frac{ab}{w}$};
\node[right]        at (C3) {$\frac{axz}{w(x+y)}$};
\node[below left]        at (D3) {$\frac{axy}{w(x+y)}$};

\node[left]        at (A4) {$x$};
\node[above left]        at (B4) {$z$};
\node[right]        at (C4) {$y$};
\node[below left]        at (D4) {$w$};

\end{tikzpicture}
\caption{A rowmotion family for $p=2$ and $q=2$.
The values $x_{i,j,k}$ are shown for all $0 \leq k \leq 4$
and all $i,j$.
The $k$-coordinate in $x_{i,j,k}$ increases in the upwards direction; the
other coordinates are oblique.}
\label{fig.2x2x5}
\end{figure}
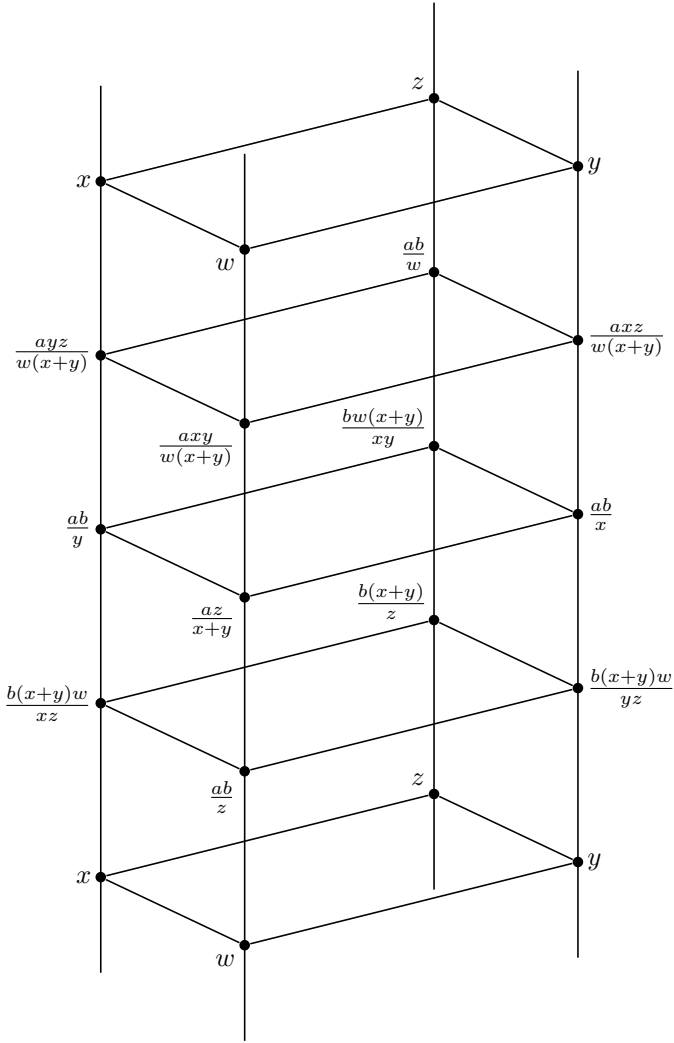

In the following decade, several proofs of these conjectures appeared, all of
which relied on detropicalization. The equation (\ref{eq.pwr}) gets replaced
by
\begin{equation}
x_{i,j,k+1}x_{i,j,k}=\left(  x_{i,j-1,k}+x_{i-1,j,k}\right)  \dfrac{1}%
{\dfrac{1}{x_{i,j+1,k+1}}+\dfrac{1}{x_{i+1,j,k+1}}}
\label{eq.brr}
\end{equation}
(with boundary conventions modified in the obvious way),
which is still uniquely solvable provided that all denominators are invertible
(as they are, e.g., when the starting values $x_{i,j,0}$ are positive reals).
(Figure~\ref{fig.2x2x5} shows five adjacent ``levels'' $k=0,1,2,3,4$ for
$p=2$ and $q=2$.)
The first proofs in \cite{BR2} recognized this new equation
\ref{eq.brr} as a relation
between certain minors of a generic $p\times\left(  p+q\right)  $-matrix,
which follows from the Pl\"{u}cker identities; a Zariski density argument then
shows that almost every solution of (\ref{eq.brr}) comes from an evaluation of
these matrix minors. A similar argument, but based on the properties of the
octahedron recurrence rather than matrix algebra, was given in
\cite{JosLiu-oct}. Explicit (if rather complicated) formulas for $x_{i,j,k}$
were proved in \cite{MusRob}, which yielded a new proof of the conjectures. In
\cite{NCBR}, the conjectures were extended to noncommutative rings (in a
modified form) and reproved using nothing but basic algebra (in the
high-school sense, albeit with summation signs and noncommutativity to spice
things up). Most of these proofs would not have worked in the original setting
of (\ref{eq.pwr}); they rely on subtraction, which was made possible by
detropicalization.\footnote{%
This actually inspired a question that is still unresolved: Do the
periodicity and reciprocity conjectures hold if the $x_{i,j,k}$ live in
a noncommutative semiring?}

\section{Outlook}

Compositions of piecewise-linear functions, such as the ones we have seen
in the previous section (and implicitly in the others as well), are not the
domain of combinatorics alone.
An application that has garnered popularity in recent years is to
\emph{feedforward neural networks} (with
\emph{ReLU activation functions}).
These compute piecewise-linear functions, i.e., rational functions over
the tropical semiring (see \cite{ZNL});
thus, these can be detropicalized as well.
However, while rowmotion is an example of such a network, its periodicity
(and general ``tameness'') appears to be rather the opposite of what is
desired of neural networks in current applications.
Nevertheless, the tropical view of neural networks has led to some
results both theoretical (\cite{BLM-TNN}) and algorithmic (\cite{LWCBM}).

Laurent phenomena are still an active field of research (see, e.g.,
\cite[\S 1.6]{BazMat}).
The same holds for birational rowmotion, which can be defined on any
finite poset (not just rectangles) and sometimes has different behavior
in the combinatorial, piecewise-linear and birational cases (see
\cite{BR1, BR2, Okada, Roby16}),
as well as several related maps such as the birational
Lalanne--Kreweras involution \cite{HopkinsJoseph} and
Coxetermotion \cite{Okada}.

\section{Appendix}

We promised a solution for the $n$ fractions puzzle. The following solution
applies to any commutative semiring $\mathbb{K}$ and any $n$ elements
$a_{1},a_{2},\ldots,a_{n}\in\mathbb{K}$ that satisfy $a_{1}a_{2}\cdots
a_{n}=1$. Assume that all denominators in (\ref{eq.nfr}) are invertible.

Set $x_{i}:=a_{1}a_{2}\cdots a_{i-1}$ for each $i\geq1$ (so that $x_{1}=1$).
From $a_{1}a_{2}\cdots a_{n}=1$, we see that the sequence $\left(  x_{1}%
,x_{2},x_{3},\ldots\right)  $ is periodic with $x_{i}=x_{i+n}$ for each $i$.
Furthermore, each $i\geq1$ and $m\geq0$ satisfy $x_{i}a_{i}a_{i+1}\cdots
a_{i+m-1}=x_{i+m}$ and thus $a_{i}a_{i+1}\cdots a_{i+m-1}=x_{i}^{-1}x_{i+m}$.
Hence, for each $i\geq1$, we have
\begin{align*}
&  1+a_{i}+a_{i}a_{i+1}+\cdots+a_{i}a_{i+1}\cdots a_{i+n-2}\\
&  =x_{i}^{-1}x_{i}+x_{i}^{-1}x_{i+1}+x_{i}^{-1}x_{i+2}+\cdots+x_{i}%
^{-1}x_{i+n-1}\\
&  =x_{i}^{-1}\left(  x_{i}+x_{i+1}+x_{i+2}+\cdots+x_{i+n-1}\right)  \\
&  =x_{i}^{-1}\left(  x_{1}+x_{2}+\cdots+x_{n}\right)
\end{align*}
by the periodicity of $\left(  x_{1},x_{2},x_{3},\ldots\right)  $. Hence,%
\begin{align*}
&  \sum_{i=1}^{n}\dfrac{1}{1+a_{i}+a_{i}a_{i+1}+\cdots+a_{i}a_{i+1}\cdots
a_{i+n-2}}\\
&  =\sum_{i=1}^{n}\dfrac{1}{x_{i}^{-1}\left(  x_{1}+x_{2}+\cdots+x_{n}\right)
}=\sum_{i=1}^{n}\dfrac{x_{i}}{x_{1}+x_{2}+\cdots+x_{n}}=1.
\end{align*}
This solves the puzzle.

Actually, the above proof works even for noncommutative semirings $\mathbb{K}%
$, under the assumption that $a_{1},a_{2},\ldots,a_{n}$ are invertible (so
that $a_{1}a_{2}\cdots a_{n}=1$ entails that $a_{i}a_{i+1}\cdots a_{i+n-1}=1$
for all $i\geq1$).

\section{Thanks}

We thank Arkadiy Berenstein,
David Einstein, Mike Joseph, Gregg Musiker, Alex Postnikov,
James Propp, and Pavlo Pylyavskyy for
enlightening and inspiring conversations throughout our tropical research.
Significant parts of that research were done and written up during our stays
at the Mathematisches Forschungsinstitut Oberwolfach in 2020 and 2023;
we thank the institute for its hospitality.

\bibliography{bibliography}

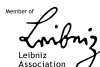
\end{document}